\documentclass[12pt]{article}
\usepackage{amsfonts}
\usepackage{amsmath}
\usepackage{graphicx}

\input tcilatex
\input tcilatex
\begin{document}

\title{}
\author{Ozlem Ersoy Hepson \\
{\scriptsize Department of Mathematics \& Computer, Eskisehir Osmangazi
University, 26480, Eskisehir, Turkey.} }
\title{\textbf{An Extended Cubic B-spline Finite Element Method for Solving
Generalized Burgers-fisher Equation}}
\maketitle

\begin{abstract}

\textbf{PACS:}
\end{abstract}

\section{Introduction}

Many numerical method uses the basis functions to solve the differential
eauations. One of the widely -used basis function are B-splines which are
used for setting up function approximation, computed aided-design and
solutions of the differential equations. There exist B-spline based
numerical methods in the numerical analysis field. \ New variants of spline
functions have been developed to form better approximation for numerical
methods. The extended B-spline functions(EBF), introduced by Han and Liu\cite%
{hanli,xu}, consist of adding additional terms to the existing B-spline
functions. Additional terms include free parameters which causes to obtain
different shapes of the B-spline form. \ EBF have the same continuity with
its standart B-spline functions and is a piecewise polynomial function of
degree 4. \ In our study we emply the collocation method to solve the
generalized Burgers-Fisher equation(GBFE). The approximation function in the
collocation method will be consist of the combination of theEBF over the
problem domain. We will observe the accuracy of the numerical solutions when
the free parameter is changed. The effect of present method will be sought \
for solutions of GBFE. Recently, EBF has started to form the numerical
methods to solve differential equations. Numerical solutions of ordinary
differential equations in the form of linear two-point boundary value
problems\cite{ex1,ex2,ex3} are given by the EBF collocation method. The
method of the collocation based on the EBF is described for solving a
one-dimensional heat equation with a nonlocal initial condition, Newell
Whitehead Segel type equation, Modified Reqularized Long Wave equation and
Advection-Diffusion equation in the works \cite{he,wk,id,di}. An extended
modified cubic B-Spline differential quadrature method is proposed to
approximate the solution of the nonlinear Burgers' equation\cite{mt}.

\noindent This study aims to construct an algorithm of combination of
Crank-Nicolson and finite element method based on extended B-spline
functions to the solutions of some initial boundary value problems defined
for the generalized Burgers-Fisher equation of the form%
\begin{equation}
u_{t}+\alpha u^{q}u_{x}-\mu u_{xx}=\eta u(1-u^{q}),\,x\in \lbrack
0,1],\,t\geq 0  \label{1}
\end{equation}%
with the initial condition%
\begin{equation}
u(x,0)=u_{0}(x)  \label{2}
\end{equation}%
and the boundary conditions%
\begin{equation}
u(0,t)=\zeta _{1}(t),\text{ }t\geq 0  \label{3}
\end{equation}%
and%
\begin{equation}
u(1,t)=\zeta _{2}(t),\text{ }t\geq 0.  \label{4}
\end{equation}%
where $\alpha ,$ $\eta $ and $q$ are parameters. The equation is a general
form form both Burgers' and Fisher equations. The solution profiles of the
equation covers many different types like single and multiple soliton
solutions.

Studies \ are going on finding both analytical and numerical solutions of
BFE. Due to including nonlinearity in the BFE, effective numerical algorithm
are necessary to understand some physical phenomena related to the BFE.
Various numerical approaches have also been implemented to the solutions of
BFE. In this paper we only mention spline related numerical scheme.The cubic
B-spline quasi-interpolation combined with low order finite difference were
suggested to solve Burgers-Fisher equation numerically. A numerical method
based on exponential spline and finite difference approximations is
developed to solve the generalized Burgers'-Fisher equation\cite{rm}. The
classical polynomial B-splines of cubic degree were used as basis to develop
a collocation method for the numerical solutions of the Burgers-Fisher
equation\cite{mit}. In \cite{gbf} generalized Burgers Fisher equation have
been solved numerically by way of the exponential cubic B-spline collocation
method.

\section{Extended Cubic B-spline Collocation Method}

The blending function of the extended cubic uniform B-spline with degree 4, $%
E_{i}(x)$, can be defined as \cite{xu}

\begin{equation}
E_{i}(x)=\frac{1}{24h^{4}}\left \{ 
\begin{array}{ll}
4h(1-\lambda )(x-x_{i-2})^{3}+3\lambda (x-x_{i-2})^{4}, & \left[
x_{i-2},x_{i-1}\right] , \\ 
\begin{array}{l}
(4-\lambda )h^{4}+12h^{3}(x-x_{i-1})+6h^{2}(2+\lambda )(x-x_{i-1})^{2} \\ 
-12h(x-x_{i-1})^{3}-3\lambda (x-x_{i-1})^{4}%
\end{array}
& \left[ x_{i-1},x_{i}\right] , \\ 
\begin{array}{l}
(4-\lambda )h^{4}-12h^{3}(x-x_{i+1})+6h^{2}(2+\lambda )(x-x_{i+1})^{2} \\ 
+12h(x-x_{i+1})^{3}-3\lambda (x-x_{i+1})^{4}%
\end{array}
& \left[ x_{i},x_{i+1}\right] , \\ 
4h(\lambda -1)(x-x_{i+2})^{3}+3\lambda (x-x_{i+2})^{4}, & \left[
x_{i+1},x_{i+2}\right] , \\ 
0 & \text{otherwise.}%
\end{array}%
\right.  \label{5}
\end{equation}

In relation (\ref{5}), the free parameter $\lambda $ is used to obtain
different form of extended cubic B-Spline functions. Note that when $\lambda
=0$, the basis function reduces to that of the cubic uniform B-spline. Graph
of the extended cubic B-splines over the interval $[0,1]$ is depicted in
Fig. 1-2 for $\lambda =-1,$ $-0.5,$ $0,$ $0.5,$ $1$ and $\lambda =-10,$ $-5,$
$0,$ $5,$ $10$ respectively.%
\begin{equation*}
\begin{array}{c}
\FRAME{itbpF}{6.2777in}{6.2777in}{0in}{}{}{fig1.jpg}{\special{language
"Scientific Word";type "GRAPHIC";maintain-aspect-ratio TRUE;display
"USEDEF";valid_file "F";width 6.2777in;height 6.2777in;depth
0in;original-width 6.25in;original-height 6.25in;cropleft "0";croptop
"1";cropright "1";cropbottom "0";filename 'Fig1.jpg';file-properties
"XNPEU";}} \\ 
\text{Fig.1: Extended cubic B-splines over the interval }[0,1]\text{ for }%
\lambda =-1,-0.5,0,0.5,1%
\end{array}%
\end{equation*}%
\begin{equation*}
\begin{array}{c}
\FRAME{itbpF}{6.2777in}{6.2777in}{0in}{}{}{fig2.jpg}{\special{language
"Scientific Word";type "GRAPHIC";display "USEDEF";valid_file "F";width
6.2777in;height 6.2777in;depth 0in;original-width 6.25in;original-height
6.25in;cropleft "0";croptop "1";cropright "1";cropbottom "0";filename
'Fig2.jpg';file-properties "XNPEU";}} \\ 
\text{Fig.2: Extended cubic B-splines over the interval }[0,1]\text{ for }%
\lambda =\lambda =-10,-5,0,5,10%
\end{array}%
\end{equation*}

$\{E_{-1}(x),E_{0}(x),\cdots ,E_{N+1}(x)\}$ forms a basis for the functions
defined over the interval $[a,b]$. Each basis function $E_{i}(x)$ is twice
continuously differentiable. The values of $E_{i}(x),E_{i}^{^{\prime }}(x),$ 
$E_{i}^{^{\prime \prime }}(x)$ at the nodal points $x_{i}$ 's computed from
Eq.(\ref{5}) are shown Table 1.

\begin{equation*}
\begin{tabular}{llllcc}
\multicolumn{6}{l}{\textbf{Table 1: }Values of $E_{i}(x)$ and its principle}
\\ 
\multicolumn{6}{l}{two derivatives at the knot points} \\ \hline
$x$ & $x_{i-2}$ & $x_{i-1}$ & $x_{i}$ & $x_{i+1}$ & $x_{i+2}$ \\ \hline
$24E_{i}$ & \multicolumn{1}{c}{$0$} & \multicolumn{1}{c}{$4-\lambda $} & 
\multicolumn{1}{c}{$16+2\lambda $} & $4-\lambda $ & $0$ \\ 
$2hE_{i}^{^{\prime }}$ & \multicolumn{1}{c}{$0$} & \multicolumn{1}{c}{$-1$}
& \multicolumn{1}{c}{$0$} & $1$ & $0$ \\ 
$2h^{2}E_{i}^{^{\prime \prime }}$ & \multicolumn{1}{c}{$0$} & 
\multicolumn{1}{c}{$2+\lambda $} & \multicolumn{1}{c}{$-4-2\lambda $} & $%
2+\lambda $ & $0$ \\ \hline
\end{tabular}%
\end{equation*}%
Suppose that the problem domain is $[a,b]$ is divided by the knots%
\begin{equation*}
\pi :a=x_{0}<x_{1}<\ldots <x_{N}=b
\end{equation*}%
into elements $[x_{i},x_{i+1}],$ $i=0,1,...,N-1$ and with mesh spacing $%
h=x_{i+1}-x_{i}=(b-a)/N,$ $i=0,1,...,N-1.$ Then the approximate solution $U$
to the unknown $u$ is written in terms of the expansion of the extended
cubic B-Spline as

\begin{equation}
U(x,t)=\sum_{i=-1}^{N+1}\delta _{i}E_{i}(x)  \label{7}
\end{equation}%
where $\delta _{i}$ are the unknown real constants and $E_{i}(x)$ are the
basis function of the extended cubic uniform B-spline. The nodal values $U$
and its first and second derivatives at the knots can be found from the (\ref%
{7}) as 
\begin{equation}
\begin{tabular}{l}
$U_{i}=U(x_{i},t)=\dfrac{4-\lambda }{24}\delta _{i-1}+\dfrac{8+\lambda }{12}%
\delta _{i}+\dfrac{4-\lambda }{24}\delta _{i+1},$ \\ 
$U_{i}^{\prime }=U^{\prime }(x_{i},t)=\dfrac{-1}{2h}\left( \delta
_{i-1}-\delta _{i+1}\right) $ \\ 
$U_{i}^{\prime \prime }=U^{\prime \prime }(x_{i},t)=\dfrac{2+\lambda }{2h^{2}%
}\left( \delta _{i-1}-2\delta _{i}+\delta _{i+1}\right) $%
\end{tabular}
\label{8}
\end{equation}

The Crank--Nicolson\ scheme is used to discretize time variables of the
unknown $U$ in the GBFE equation so that one obtain the time discretized
form of the equation as 
\begin{equation}
\frac{U^{n+1}-U^{n}}{\Delta t}=-\alpha \frac{%
(U^{q}U_{x})^{n+1}+(U^{q}U_{x})^{n}}{2}+\mu \frac{U_{xx}^{n+1}+U_{xx}^{n}}{2}%
+\eta \frac{U^{n+1}+U^{n}}{2}-\eta \frac{\left( U^{q+1}\right) ^{n+1}+\left(
U^{q+1}\right) ^{n}}{2}  \label{9}
\end{equation}%
where $U^{n+1}=U(x,t^{n+1})$ is the solution of the equation at the $(n+1)$%
th. time level. \ Here $t^{n+1}$ $=t^{n}+\Delta t$ and $\Delta t$ is the
time step, superscripts denote $n$ th time level , $t^{n}=n\Delta t.$

The nonlinear term $(U^{q}U_{x})^{n+1}$ and $\left( U^{q+1}\right) ^{n+1}$
in Eq. (\ref{9}) is linearized by using the following form \cite{rubin}:%
\begin{equation*}
(U^{q}Ux)^{n+1}=\left( U^{q}\right) ^{n}U_{x}^{n+1}+q\left( U^{q-1}\right)
^{n}U_{x}^{n}U^{n+1}-q\left( U^{q}\right) ^{n}U_{x}^{n}
\end{equation*}%
and%
\begin{eqnarray*}
(U^{q+1})^{n+1} &=&(U^{q}U)^{n+1} \\
&=&\left( U^{q}\right) ^{n}U^{n+1}+q\left( U^{q-1}\right)
^{n}U^{n}U^{n+1}-q\left( U^{q}\right) ^{n}U^{n} \\
&=&(1+q)\left( U^{q}\right) ^{n}U^{n+1}-q\left( U^{q+1}\right) ^{n}
\end{eqnarray*}%
So Eq. (\ref{9}) is discretized in time as%
\begin{eqnarray}
&&U^{n+1}+\alpha \frac{\Delta t}{2}(L_{1})^{q}U_{x}^{n+1}+\alpha \frac{%
\Delta t}{2}q(L_{1})^{q-1}L_{2}U^{n+1}-\mu \frac{\Delta t}{2}%
U_{xx}^{n+1}-\eta \frac{\Delta t}{2}U^{n+1}  \notag \\
&&+\eta \frac{\Delta t}{2}(1+q)(L_{1})^{q}U^{n+1}  \notag \\
&=&U^{n}-\alpha \frac{\Delta t}{2}(1-q)(L_{1})^{q}U_{x}^{n}+\mu \frac{\Delta
t}{2}U_{xx}^{n}+\eta \frac{\Delta t}{2}U^{n}-\eta \frac{\Delta t}{2}%
(1-q)(L_{1})^{q}U^{n}  \label{10}
\end{eqnarray}%
Substitution (\ref{8}) into (\ref{10}) leads to the fully-discretized
equation:

\begin{equation}
\begin{array}{l}
\left[ \left( 1+\alpha \dfrac{\Delta t}{2}q(L_{1})^{q-1}L_{2}-\eta \dfrac{%
\Delta t}{2}+\eta \dfrac{\Delta t}{2}(1+q)(L_{1})^{q}\right) \alpha
_{1}+\alpha \dfrac{\Delta t}{2}(L_{1})^{q}\beta _{1}-\mu \dfrac{\Delta t}{2}%
\gamma _{1}\right] \delta _{i-1}^{n+1} \\ 
+\left[ \left( 1+\alpha \dfrac{\Delta t}{2}q(L_{1})^{q-1}L_{2}-\eta \dfrac{%
\Delta t}{2}+\eta \dfrac{\Delta t}{2}(1+q)(L_{1})^{q}\right) \alpha _{2}-\mu 
\dfrac{\Delta t}{2}\gamma _{2}\right] \delta _{i}^{n+1} \\ 
+\left[ \left( 1+\alpha \dfrac{\Delta t}{2}q(L_{1})^{q-1}L_{2}-\eta \dfrac{%
\Delta t}{2}+\eta \dfrac{\Delta t}{2}(1+q)(L_{1})^{q}\right) \alpha
_{1}-\alpha \dfrac{\Delta t}{2}(L_{1})^{q}\beta _{1}-\mu \dfrac{\Delta t}{2}%
\gamma _{1}\right] \delta _{i+1}^{n+1} \\ 
=\left[ \left( 1+\eta \dfrac{\Delta t}{2}-\eta \dfrac{\Delta t}{2}%
(1-q)(L_{1})^{q}\right) \alpha _{1}-\alpha \dfrac{\Delta t}{2}%
(1-q)(L_{1})^{q}\beta _{1}+\mu \dfrac{\Delta t}{2}\gamma _{1}\right] \delta
_{i-1}^{n} \\ 
+\left[ \left( 1+\eta \dfrac{\Delta t}{2}-\eta \dfrac{\Delta t}{2}%
(1-q)(L_{1})^{q}\right) \alpha _{2}+\mu \dfrac{\Delta t}{2}\gamma _{2}\right]
\delta _{i}^{n} \\ 
+\left[ \left( 1+\eta \dfrac{\Delta t}{2}-\eta \dfrac{\Delta t}{2}%
(1-q)(L_{1})^{q}\right) \alpha _{1}+\alpha \dfrac{\Delta t}{2}%
(1-q)(L_{1})^{q}\beta _{1}+\mu \dfrac{\Delta t}{2}\gamma _{1}\right] \delta
_{i+1}^{n}%
\end{array}
\label{12}
\end{equation}%
where%
\begin{eqnarray*}
L_{1} &=&\alpha _{1}\delta _{i-1}^{n}+\alpha _{2}\delta _{i}^{n}+\alpha
_{1}\delta _{i+1}^{n} \\
L_{2} &=&\beta _{1}\delta _{i-1}^{n}-\beta _{1}\delta _{i+1}^{n}
\end{eqnarray*}%
\begin{eqnarray*}
\alpha _{1} &=&\dfrac{4-\lambda }{24},\text{ }\alpha _{2}=\dfrac{8+\lambda }{%
12} \\
\beta _{1} &=&-\dfrac{1}{2h},\text{ }\gamma _{1}=\dfrac{2+\lambda }{2h^{2}},%
\text{ }\gamma _{2}=-\dfrac{4+2\lambda }{2h^{2}}
\end{eqnarray*}

A linear system of $(N+1)$ equations in $N+3$ unknown is obtained. However,
two additional linear equations are needed to obtain the values of $n+3$
unknown parameters $\mathbf{d}^{n+1}=(\delta _{-1}^{n+1},\delta
_{0}^{n+1},\ldots ,\delta _{N+1}^{n+1})$. \ To make solvable the system,
boundary conditions $U_{0}=\zeta _{1},$ $U_{N}=\zeta _{2}$ are used to find
two additional linear equations: 
\begin{eqnarray}
\delta _{-1} &=&\frac{1}{\alpha _{1}}\left( U_{0}-\alpha _{2}\delta
_{0}-\alpha _{3}\delta _{1}\right) ,  \label{13} \\
\delta _{N+1} &=&\frac{1}{\alpha _{3}}\left( U_{n}-\alpha _{1}\delta
_{N-1}-\alpha _{2}\delta _{N}\right) .  \notag
\end{eqnarray}%
(\ref{13}) can be used to eliminate $\delta _{-1},$ $\delta _{N+1}$ from the
system (\ref{12}) which then becomes the solvable matrix equation. A variant
of Thomas algorithm is used to solve the system.

Before starting the iteration of the Eq. (\ref{12}), initial parameters $%
\delta _{-1}^{0},\delta _{0}^{0},\ldots ,\delta _{N+1}^{0}$ must be
determined from the initial condition and first space derivative of the
initial conditions at the boundaries as the following:

\begin{enumerate}
\item $U(x_{i},0)$ $=U(x_{i},0),$ $i=0,...,N$

\item $(U_{x})(x_{0},0)=U(x_{0})$

\item $(U_{x})(x_{N},0)=U(x_{N}).$
\end{enumerate}

\section{Numerical tests}

In this section, some numerical solutions of the GBFE with the extended
cubic B-Spline collocation are presented. To show the efficiency of the
present method for our problem in comparison with the exact solution we
report maximum error which is defined by

\begin{equation*}
L_{\infty }=\left \vert u-U\right \vert _{\infty }=\max \limits_{j}\left
\vert u_{j}-U_{j}^{n}\right \vert
\end{equation*}%
where $U$ is the solution obtained by Eq. (\ref{1}), solved by the extended
cubic B-spline collocation method and $u$ is the exact solution.

\subsection{Example 1}

Consider GBFE with the initial condition%
\begin{equation*}
u(x,0)=\left \{ \frac{1}{2}+\frac{1}{2}\tanh (\frac{-\alpha q}{2(q+1)}%
)x\right \} ^{\dfrac{1}{q}}=\varphi (x)
\end{equation*}%
the boundary conditions%
\begin{equation*}
u(0,t)=\left( \frac{1}{2}+\frac{1}{2}\tanh \left[ \frac{-\alpha q}{2(q+1)}%
\left( -\left( \frac{\alpha }{q+1}+\frac{\eta (q+1)}{\alpha }\right)
t\right) \right] \right) ^{\dfrac{1}{q}}=\zeta _{1}(t),\text{ }t\geq 0
\end{equation*}%
and%
\begin{equation*}
u(1,t)=\left( \frac{1}{2}+\frac{1}{2}\tanh \left[ \frac{-\alpha q}{2(q+1)}%
\left( 1-\left( \frac{\alpha }{q+1}+\frac{\eta (q+1)}{\alpha }\right)
t\right) \right] \right) ^{\dfrac{1}{q}}=\zeta _{2}(t),\text{ }t\geq 0
\end{equation*}%
its exact solution%
\begin{equation*}
u(x,t)=\left( \frac{1}{2}+\frac{1}{2}\tanh \left[ \frac{-\alpha q}{2(q+1)}%
\left( x-\left( \frac{\alpha }{q+1}+\frac{\eta (q+1)}{\alpha }\right)
t\right) \right] \right) ^{\dfrac{1}{q}},\text{ }t\geq 0
\end{equation*}
We run the program for three sets of parameters to make comparison with
results of some earlier studies \cite{bf,gold,ZZ,BF1} .

We show discrete $L_{\infty }$ error norms for$,\alpha =0.1,$ $\eta
=-0.0025, $ $\Delta t=0.0001,N=16$ and results are documented in Table 2.
Absolute error of solution at $t=1$ is depicted in Fig. 3 when $q=1,$

\begin{equation*}
\begin{tabular}{|l|}
\hline
Table 2: $L_{\infty }$ error norms \\ \hline
$%
\begin{tabular}{lllll}
& Time($t$) & $l=0$ & $(l=-0.000003)$ & \cite{bf} \\ \hline
$q=1$ & $0.1$ & $1.08646\times 10^{-12}$ & $1.02251\times 10^{-13}$ & $%
1.32396\times 10^{-11}$ \\ 
& $0.2$ & $1.46944\times 10^{-12}$ & $1.24456\times 10^{-13}$ & $%
1.78026\times 10^{-11}$ \\ 
& $0.3$ & $1.61926\times 10^{-12}$ & $1.24456\times 10^{-13}$ & $%
1.94258\times 10^{-11}$ \\ 
& $0.4$ & $1.67277\times 10^{-12}$ & $1.24456\times 10^{-13}$ & $%
2.00083\times 10^{-11}$ \\ 
& $0.5$ & $1.67277\times 10^{-12}$ & $1.24456\times 10^{-13}$ & $%
2.02158\times 10^{-11}$ \\ 
&  &  &  &  \\ \hline
$q=2$ & $0.1$ & $2.17542\times 10^{-11}$ & $1.34003\times 10^{-13}$ & $%
2.84700\times 10^{-10}$ \\ 
& $0.2$ & $3.02457\times 10^{-11}$ & $1.47881\times 10^{-13}$ & $%
3.87950\times 10^{-10}$ \\ 
& $0.3$ & $3.33861\times 10^{-11}$ & $1.47881\times 10^{-13}$ & $%
4.24646\times 10^{-10}$ \\ 
& $0.4$ & $3.45414\times 10^{-11}$ & $1.47104\times 10^{-13}$ & $%
4.37589\times 10^{-10}$ \\ 
& $0.5$ & $3.49593\times 10^{-11}$ & $1.43551\times 10^{-13}$ & $%
4.02050\times 10^{-10}$ \\ 
&  &  &  &  \\ 
$q=4$ & $0.1$ & $3.12324\times 10^{-11}$ & $2.65165\times 10^{-12}$ & $%
3.99168\times 10^{-10}$ \\ 
& $0.2$ & $4.34227\times 10^{-11}$ & $3.67783\times 10^{-12}$ & $%
5.43802\times 10^{-10}$ \\ 
& $0.3$ & $4.79300\times 10^{-11}$ & $4.05153\times 10^{-12}$ & $%
5.95169\times 10^{-10}$ \\ 
& $0.4$ & $4.95853\times 10^{-11}$ & $4.18043\times 10^{-12}$ & $%
6.13233\times 10^{-10}$ \\ 
& $0.5$ & $5.01746\times 10^{-11}$ & $4.21829\times 10^{-12}$ & $%
6.19407\times 10^{-10}$%
\end{tabular}%
$ \\ \hline
\end{tabular}%
\end{equation*}

\begin{equation*}
\begin{array}{c}
\\ 
\text{Figure 3: The absolute errors }\alpha =0.1,\text{ }\eta =-0.0025,\text{
}\Delta t=0.0001,\text{ }q=1%
\end{array}%
\end{equation*}

For some values of $q$ and $t$, discrete $L_{\infty }$ error norms are
recorded for parametes $\alpha =1,$ $\eta =1,$ $\Delta t=0.0001$ in table 3$%
. $ Absolute error for $q=1$ at $t=1$ is drawn in Figure 4.

\begin{equation*}
\begin{tabular}{|l|}
\hline
Table 3: $L_{\infty }$ error norms \\ \hline
$%
\begin{tabular}{lllll}
& Time($t$) & $l=0$ & Various $l$ & \cite{bf} \\ 
$q=1$ & $0.2$ & $5.58038\times 10^{-8}$ & $1.94765\times
10^{-10}(p=-0.000319)$ & $5.55746\times 10^{-7}$ \\ 
& $0.4$ & $8.54479\times 10^{-8}$ & $1.54361\times 10^{-9}$ & $9.05507\times
10^{-7}$ \\ 
& $0.6$ & $2.04362\times 10^{-7}$ & $1.37196\times 10^{-8}$ & $2.18808\times
10^{-6}$ \\ 
& $0.8$ & $2.80869\times 10^{-7}$ & $4.72669\times 10^{-8}$ & $2.93314\times
10^{-7}$ \\ 
& $1.0$ & $2.99588\times 10^{-7}$ & $1.02753\times 10^{-7}$ & $3.01455\times
10^{-6}$ \\ 
&  &  &  &  \\ 
$q=2$ & $0.2$ & $2.82618\times 10^{-7}$ & $7.64068\times
10^{-9}(p=-0.0002930)$ & $2.56108\times 10^{-6}$ \\ 
& $0.4$ & $4.62302\times 10^{-7}$ & $2.75394\times 10^{-8}$ & $4.24308\times
10^{-6}$ \\ 
& $0.6$ & $4.29008\times 10^{-7}$ & $7.20780\times 10^{-8}$ & $3.56848\times
10^{-6}$ \\ 
& $0.8$ & $2.56283\times 10^{-7}$ & $2.07798\times 10^{-7}$ & $1.46518\times
10^{-6}$ \\ 
& $1.0$ & $8.03168\times 10^{-8}$ & $3.03531\times 10^{-7}$ & $5.54230\times
10^{-6}$ \\ 
&  &  &  &  \\ 
$q=4$ & $0.2$ & $3.98349\times 10^{-7}$ & $6.38513\times 10^{-7}(p=0.000193)$
& $1.76161\times 10^{-6}$ \\ 
& $0.4$ & $2.64952\times 10^{-7}$ & $5.11063\times 10^{-7}$ & $4.17351\times
10^{-7}$ \\ 
& $0.6$ & $1.73948\times 10^{-8}$ & $1.59753\times 10^{-7}$ & $2.42401\times
10^{-6}$ \\ 
& $0.8$ & $8.61362\times 10^{-8}$ & $3.38624\times 10^{-9}$ & $2.35757\times
10^{-6}$ \\ 
& $1.0$ & $6.63329\times 10^{-7}$ & $2.41389\times 10^{-8}$ & $1.44350\times
10^{-6}$%
\end{tabular}%
$ \\ \hline
\end{tabular}%
\end{equation*}%
\begin{equation*}
\begin{array}{c}
\\ 
\text{Fig. 3: The absolute errors }\alpha =1,\text{ }\eta =1,\text{ }\Delta
t=0.0001,\text{ }q=1,\text{ }t=1,l=0%
\end{array}%
\end{equation*}

Last, the program is rerun with different parameters $\alpha =0.01,$ $0.001,$
$\eta =1,10,100$ time step $\Delta t=0.0001,$ $q=1,$space step $N=8$ to make
comparison with results of the spectral collocation method the discontinues
Galerkin method and cubic B-spline collocation documented in Table 5 at
times $t=0.5.$%
\begin{equation*}
\begin{tabular}{|l|}
\hline
Table 4: $L_{\infty }$ error norm for the solutions of Example 1 at $t=0.5$
for $\alpha =0.001$ and $\alpha =0.0001$ \\ \hline
$%
\begin{tabular}{lllllll}
\hline
& $\eta $ & Present($l=0$) & Various $l$ & \cite{gold} & \cite{ZZ} & \cite%
{BF1} \\ \hline
$\alpha =0.01$ & $1$ & $2.43664\times 10^{-11}$ & $2.95541\times
10^{-13}(p=0.00034)$ & $4.6763\times 10^{-12}$ & $2.8999\times 10^{-13}$ & $%
2.4264\times 10^{-12}$ \\ 
& $10$ & $6.89380\times 10^{-10}$ & $1.18782\times 10^{-12}(p=0.07066)$ & $%
6.2529\times 10^{-12}$ & $3.3184\times 10^{-13}$ & $1.2833\times 10^{-13}$
\\ 
& $100$ & $9.32587\times 10^{-15}$ & $4.44089\times 10^{-16}(p=-0.56755)$ & $%
8.0269\times 10^{-12}$ & $2.4225\times 10^{-13}$ & $1.2500\times 10^{-12}$
\\ 
&  &  &  &  &  &  \\ \hline
$\alpha =0.001$ & $1$ & $2.43607\times 10^{-11}$ & $4.138911\times
10^{-13}(p=0.03311)$ & $4.5374\times 10^{-12}$ & $2.8821\times 10^{-13}$ & $%
2.4251\times 10^{-12}$ \\ 
& $10$ & $6.87854\times 10^{-10}$ &  & $6.0540\times 10^{-12}$ & $%
3.3295\times 10^{-13}$ & $1.2832\times 10^{-13}$ \\ 
& $100$ & $9.10383\times 10^{-15}$ & $8.88178\times 10^{-16}(p=-0.85953)$ & $%
8.1424\times 10^{-13}$ & $2.4480\times 10^{-13}$ & $1.2500\times 10^{-12}$%
\end{tabular}%
$ \\ \hline
\end{tabular}%
\end{equation*}

\subsection{Example 2}

With initial profile as 
\begin{equation*}
u(x,0)=\exp (-40x^{2})
\end{equation*}%
evolution of solution of GBFE is depicted with different parameters $\alpha ,
$ $\eta $ and $\mu $ for $N=80,$ $\Delta t=0.001,$ $t=1.5$ in Figure 4-7%
\begin{equation*}
\begin{tabular}{ll}
\FRAME{itbpF}{3.7075in}{3.7075in}{0in}{}{}{rc1.bmp}{\special{language
"Scientific Word";type "GRAPHIC";maintain-aspect-ratio TRUE;display
"USEDEF";valid_file "F";width 3.7075in;height 3.7075in;depth
0in;original-width 3.6599in;original-height 3.6599in;cropleft "0";croptop
"1";cropright "1";cropbottom "0";filename 'RC1.bmp';file-properties "XNPEU";}%
} & \FRAME{itbpF}{3.7075in}{3.7075in}{0in}{}{}{rc2.bmp}{\special{language
"Scientific Word";type "GRAPHIC";maintain-aspect-ratio TRUE;display
"USEDEF";valid_file "F";width 3.7075in;height 3.7075in;depth
0in;original-width 3.6599in;original-height 3.6599in;cropleft "0";croptop
"1";cropright "1";cropbottom "0";filename 'RC2.bmp';file-properties "XNPEU";}%
} \\ 
Figure 4: $\alpha =0,$ $\eta =1$ and $\mu =0.1$ & Figure 5: $\alpha =1,$ $%
\eta =0.02$ and $\mu =0.02$%
\end{tabular}%
\end{equation*}%
\begin{equation*}
\begin{tabular}{ll}
\FRAME{itbpF}{3.7075in}{3.7075in}{0in}{}{}{rc3.bmp}{\special{language
"Scientific Word";type "GRAPHIC";maintain-aspect-ratio TRUE;display
"USEDEF";valid_file "F";width 3.7075in;height 3.7075in;depth
0in;original-width 3.6599in;original-height 3.6599in;cropleft "0";croptop
"1";cropright "1";cropbottom "0";filename 'RC3.bmp';file-properties "XNPEU";}%
} & \FRAME{itbpF}{3.7075in}{3.7075in}{0in}{}{}{rc4.bmp}{\special{language
"Scientific Word";type "GRAPHIC";maintain-aspect-ratio TRUE;display
"USEDEF";valid_file "F";width 3.7075in;height 3.7075in;depth
0in;original-width 3.6599in;original-height 3.6599in;cropleft "0";croptop
"1";cropright "1";cropbottom "0";filename 'RC4.bmp';file-properties "XNPEU";}%
} \\ 
Figure 6: $\alpha =1,$ $\eta =0.02$ and $\mu =0.002$ & Figure 7: $\alpha =1,$
$\eta =0.02$ and $\mu =0.0002$%
\end{tabular}%
\end{equation*}

\subsection{Example 3}

Consider GBFE with $\eta =0,q=1$ and $\alpha =1$ with initial and boundary
conditions%
\begin{equation*}
u(x,0)=x(1-x^{2}),\text{ }0<x<1
\end{equation*}%
and%
\begin{equation*}
u(0,t)=u(1,t)=0,\text{ }t\geq 0.
\end{equation*}%
In this example, we computed solutions for $\mu =2^{-2}$ and $\mu =2^{-6}$
at $t=0.1,$ $0.3,$ $0.6,$ $0.9$ with step size $\Delta t=0.001.$ The
obtained solutions are plotted in Fig. 8-9. Next, we obtained resulst for $%
\mu =2^{-2},$ $\mu =2^{-4},$ $\mu =2^{-6},$ $\mu =2^{-8}$ and $t=0.5,$ $0.9,$
respectively, with step size $\Delta t=0.001.$ We plotted corresponding
solution curves in Fig. 10-11.%
\begin{equation*}
\begin{array}{c}
\begin{tabular}{ll}
\FRAME{itbpF}{3.3927in}{2.7095in}{0in}{}{}{3a.bmp}{\special{language
"Scientific Word";type "GRAPHIC";maintain-aspect-ratio TRUE;display
"USEDEF";valid_file "F";width 3.3927in;height 2.7095in;depth
0in;original-width 3.346in;original-height 2.6671in;cropleft "0";croptop
"1";cropright "1";cropbottom "0";filename '3a.bmp';file-properties "XNPEU";}}
& \FRAME{itbpF}{3.3927in}{2.7095in}{0in}{}{}{3b.bmp}{\special{language
"Scientific Word";type "GRAPHIC";maintain-aspect-ratio TRUE;display
"USEDEF";valid_file "F";width 3.3927in;height 2.7095in;depth
0in;original-width 3.346in;original-height 2.6671in;cropleft "0";croptop
"1";cropright "1";cropbottom "0";filename '3b.bmp';file-properties "XNPEU";}}
\\ 
Figure 8: $\alpha =1,$ $\eta =0$ and $\mu =2^{-2}$. & Figure 9: $\alpha =1,$ 
$\eta =0$ and $\mu =2^{-6}$%
\end{tabular}
\\ 
\text{Computed solutions of Example 3 for different time levels at fixed.}%
\mu \text{.}%
\end{array}%
\end{equation*}

\begin{equation*}
\begin{array}{c}
\begin{tabular}{ll}
\FRAME{itbpF}{3.3797in}{2.7095in}{0in}{}{}{3c.bmp}{\special{language
"Scientific Word";type "GRAPHIC";maintain-aspect-ratio TRUE;display
"USEDEF";valid_file "F";width 3.3797in;height 2.7095in;depth
0in;original-width 3.333in;original-height 2.6671in;cropleft "0";croptop
"1";cropright "1";cropbottom "0";filename '3c.bmp';file-properties "XNPEU";}}
& \FRAME{itbpF}{3.3866in}{2.7294in}{0in}{}{}{3d.bmp}{\special{language
"Scientific Word";type "GRAPHIC";maintain-aspect-ratio TRUE;display
"USEDEF";valid_file "F";width 3.3866in;height 2.7294in;depth
0in;original-width 3.3399in;original-height 2.687in;cropleft "0";croptop
"1";cropright "1";cropbottom "0";filename '3d.bmp';file-properties "XNPEU";}}
\\ 
Figure 10: $\alpha =1,$ $\eta =0$ and $t=0.5$. & Figure 11: $\alpha =1,$ $%
\eta =0$ and $t=0.9$%
\end{tabular}
\\ 
\text{Computed solutions of Example 3 for different time levels at fixed.}t%
\text{.}%
\end{array}%
\end{equation*}

\section{Conclusion}

The extended B-spline collocation method in space and the Crank-Nicolson
method in time were implemented for obtaining explicit solution of the
generalized Burgers--Fisher equation. An alternatif algorithm is suggested
to solve the GBFE. For first text problem, present method provided lower
error than the cubic B-spline quasi-interpolation method. Similar \ results
of the present method are obtained when compared with A spectral domain
decomposition approach, local discontinuous Galerkin method and cubic
B-spline collocation method when we catch a suitable free parameter which is
determined by scanning the predeterming interval with an small space
increment, we have achived less error mostly

\end{document}